\documentclass[twoside,12pt]{article}

\usepackage{amssymb,amsmath}

\newenvironment{proof}{\begin{trivlist}\item[]{\it
Proof.}}{\hfill$\square$\end{trivlist}}

\newtheorem{theorem}{Theorem}[section]
\newtheorem{corollary}[theorem]{Corollary}
\newtheorem{definition}[theorem]{Definition}
\newtheorem{lemma}[theorem]{Lemma}
\newtheorem{proposition}[theorem]{Proposition}
\newtheorem{remark}[theorem]{Remark}

\def\ca{{\cal A}}
\def\cb{{\cal B}}
\def\ch{{\cal O}}
\def\cc{{\cal C}}

\def\car{\ca(R)}

\def\a{\alpha} 
\def\b{\beta}

\def\l{\lambda}
\def\r{\rho}

\def\bigd{\Delta}
\def\aug{\varepsilon} 
\def\suma{\sum_{(a)}\,}
\def\sumb{\sum_{(b)}\,}
\def\sumab{\sum_{(ab)}\,}
\def\sumv{\sum_{(v)}\,}

\def\aone{a_1}
\def\atwo{a_2}
\def\athree{a_3} 
\def\bone{b_1}
\def\btwo{b_2}
\def\bthree{b_3}
\def\cone{c_1}
\def\ctwo{c_2}

\def\br{{\bf r}}
\def\brbar{{\bar{{\bf r}}}}

\def\id{{\mathrm{id}}}
\def\im{{\mathrm{im}}}

\def\mc{{\mathbb C}}
\def\mq{{\mathbb Q}}
\def\mz{{\mathbb Z}}
\def\mn{{\mathbb N}}
\def\si{\sigma_i}
\def\ti{\tau_i}

\def\slN{{\mathfrak{sl}}_N}
\def\soN{{\mathfrak{so}}_N}
\def\spN{{\mathfrak{sp}}_N}
\def\oslqn{\ch(SL_q(N))}
\def\ooqn{\ch(O_q(N))}
\def\ospqn{\ch(Sp_q(N))}
\def\osoqn{\ch(SO_q(N))}
\def\oglqn{\ch(GL_q(N))}
\def\omqn{\ch(M_q(N))}

\def\ok{k[z,z^{-1}]}
\def\ozeta{k[\zeta\mid \zeta^N=1]}

\def\detq{\mathrm{det}_q}

\def\cor{\rho\mbox{{\small -co-}}\ch}
\def\col{\lambda\mbox{{\small -co-}}\ch}
\def\coa{\alpha\mbox{{\small -co-}}\ch}
\def\cob{\beta\mbox{{\small -co-}}\ch}

\def\sse{\subseteq} 
\def\ot{\otimes} 
\def\goesto{\longrightarrow}

\date{}
\begin{document}
\title{Weakly multiplicative coactions of
\\ quantized function algebras} 
\author{M Domokos\thanks{This research was supported through a European
Community
Marie Curie Fellowship. Partially supported by OTKA No. F 32325 and
T 34530.}\; and T H Lenagan} 
\maketitle

\begin{abstract}
A condition is identified which guarantees that the coinvariants of a
coaction of a Hopf algebra on an algebra form a subalgebra,
even though the coaction may fail to be an algebra homomorphism. A
Hilbert Theorem (finite generation of the subalgebra of coinvariants) is
obtained for such coactions of a cosemisimple Hopf algebra. 
This is applied for two coactions 
$\alpha,\beta:\ca\to\ca\otimes\ch$, 
where $\ca$ is the coordinate algebra of the quantum matrix space 
associated with the quantized coordinate algebra $\ch$ of a classical group, 
and $\alpha$, $\beta$ are quantum analogues of the conjugation action on 
matrices. Provided that $\ch$ is cosemisimple and coquasitriangular, 
the $\alpha$-coinvariants and the 
$\beta$-coinvariants form two finitely generated, 
commutative, graded subalgebras of $\ca$, having the same Hilbert series. 
Consequently, the cocommutative elements and the $S^2$-cocommutative
elements in $\ch$ form finitely generated subalgebras. 
A Hopf algebra monomorphism from the quantum general linear group to
Laurent polynomials over the quantum special linear group is found and
used to explain the strong relationship between the corepresentation
(and coinvariant) theories of these quantum groups. 
\end{abstract}

\medskip
\noindent 2000 Mathematics Subject Classification. 16W30, 16W35, 20G42,
17B37, 57T05, 81R50

\noindent Keywords: coinvariant, cosemisimple Hopf algebra, quantum
general linear group, quantum special linear group, quantum orthogonal group, 
quantum symplectic group, 
coquasitriangular Hopf algebra, quantized function algebras, 
adjoint coaction, cocommutative element 

\bigskip 

Let $k$ be a field.  Let  $A$ be a $k$-algebra
and $\ch$ a Hopf algebra over $k$ and suppose that 
$\beta: A \goesto A\ot \ch$ is a right coaction.  
The {\em coinvariants} of $\beta$ consist of the set
$A^{\cob}:= \{c\in A \mid \beta(c) = c\ot 1\}$.  When $\beta$ is also an
algebra homomorphism, $A^{\cob}$ is automatically a subalgebra of $A$. 
However, it has recently become apparent that
progress can be made even when $\beta$ is not an algebra homomorphism,
see \cite{dl, dfl, gl, glr}, for example, where examples of such
coactions are studied, motivated by seeking quantum versions of results
concerning the classical invariant theory of the general linear and
special linear groups.  In these quantum cases, at the
outset, it is not even clear that the set of coinvariants forms a
subalgebra of $A$.  However, under weaker conditions it is sometimes
possible to show that the set of coinvariants does indeed form a
subalgebra, and that this subalgebra enjoys desirable properties.  
The first purpose of 
this paper is to isolate the key points that make this process
work.  In particular, we are able to show that, under certain
assumptions on $\beta$, when $\ch$ is cosemisimple, the set of
coinvariants forms a subalgebra which is finitely generated.  

Our leading examples for such coactions are associated to 
quantized algebras of functions on simple Lie groups. 
Let $\ch$ be any of these quantum versions of classical groups 
(cf. \cite{rtf}). 
We consider two right coactions $\alpha$, $\beta$ of $\ch$ on 
the coordinate algebra $\ca$ of the corresponding quantum matrix space. 
(In the case of $\oslqn$, the coordinate algebra of the 
quantum special linear group, both of these coactions are possible 
quantum deformations of the conjugation action of $SL(N)$ on 
$N\times N$ matrices). As an application of the general 
considerations mentioned above, we show that provided that 
$\ch$ is cosemisimple and coquasitriangular, 
each of the 
$\alpha$-coinvariants and $\beta$-coinvariants forms a finitely generated, 
graded, commutative subalgebra of $\ca$, 
and these subalgebras have the same Hilbert series. 
In particular, both the cocommutative elements and the 
$S^2$-cocommutative elements (cf. \cite{cw}) form a finitely 
generated subalgebra in $\ch$. 

The last section clarifies the relation between $\oslqn$ and 
$\oglqn$, the coordinate algebra of the quantum general linear group, 
which has consequences for coinvariant theory. 
Levasseur and Stafford, \cite{ls}, have shown that 
$\oglqn$ and 
$\oslqn\otimes\ok$ are isomorphic as $k$-algebras, and that  
ring theoretic properties of $\oslqn$ can be derived from this
isomorphism. However, the algebra isomorphism 
$\oglqn\to\oslqn\otimes\ok$ 
given in \cite{ls} is not a morphism of coalgebras. 
Here we construct 
a Hopf algebra homomorphism which embeds $\oglqn$ into the Hopf algebra
$\oslqn\otimes\ok$ in such a way that $\oslqn\otimes\ok$ is a finitely
generated free module over the image of $\oglqn$. This explains the
strong relationship between the corepresentation theories of these
quantum groups. 
To illustrate this, we derive explicitly generating 
$\oslqn$-coinvariants in $\omqn$ and $\oslqn$  
(with respect to the coactions $\alpha$, $\beta$, 
under the assumption that $q$ is not a root of unity)  
from the corresponding results for $\oglqn$ obtained
in \cite{dl}.

\section{Weakly multiplicative coactions}\label{coinvs}
\marginpar{coinvs}

First, we identify a condition that is frequently satisfied by 
coactions even when they are not algebra maps, and which guarantees that
the coinvariants form a subalgebra. 

Let $A$ be a $k$-algebra and let $\ch$ be a Hopf algebra over $k$.  Let
$\bigd, \; S$ and $\aug$ denote the comultiplication, antipode  
and counit maps of $\ch$,
respectively.  If
$\b:A\goesto A\ot\ch$ is a right coaction, then we say that
$\b$ is {\em left weakly multiplicative} if $\b(ab) = \b(a)\b(b)$ for
all $b\in A$, provided that $\b(a) = a \ot 1$. 
Similarly, $\b$ is {\em right weakly multiplicative} if $\b(ab) =
\b(a)\b(b)$ for
all $b\in A$, provided that $\b(b) = b \ot 1$.
The following trivial 
lemma shows that it is useful to identify when a right coaction
is left or right weakly multiplicative.

\begin{lemma}\label{subalg}
Let $A$ be a $k$-algebra and let $\ch$ be a Hopf algebra over $k$.  If 
$\b:A\goesto A\ot\ch$ is either  a left or right 
weakly multiplicative right 
coaction then the set of coinvariants $A^{\cob} = \{c\in A \mid \b(c) =
c\ot 1\}$ forms a subalgebra of $A$.
\end{lemma}

\begin{proof}
All that needs to be proved is that the product of two coinvariants are a
coinvariant.  Let $a,b $ be coinvariants, so that $\b(a) = a\ot 1$ and
$\b(b) = b\ot 1$.  Then $\b(ab) = \b(a)\b(b) = (a\ot 1)(b\ot 1) = ab \ot
1$, as required, by using the left or right
weakly multiplicative condition, as appropriate. 
\end{proof}

Next, we show how to construct some left or right 
weakly multiplicative right
coactions which, in general, are not algebra homomorphisms. 

We fix the following notation. 
Let $A$ be  a $k$-algebra and let $\ch$ be a Hopf algebra over
$k$.  Suppose that  $\l: A \goesto \ch \ot A$ is a left coaction and
that $\r:A\goesto A\ot\ch$ is a right coaction, and that both $\l$ and
$\r$ 
are algebra homomorphisms. Suppose further that $\l$ and
$\r$ commute, in the sense that $(\l\ot\id)\circ\r = (\id\ot\r)\circ\l$.

Define 
\begin{equation}\label{eq:alpha}
\a :=  (\id\ot m)\circ
\tau_{132}\circ(S\ot\id\ot\id)\circ(\l\ot\id)\circ\r:A\goesto
A\ot\ch,
\end{equation} 
and 
\begin{equation}\label{eq:beta}
\b:=  (\id\ot m)\circ
(\tau\ot\id)\circ(S\ot\id\ot\id)\circ(\l\ot\id)\circ\r:A\goesto
A\ot\ch 
\end{equation} 
where $S$ is the antipode of $\ca$, while 
$\tau_{132}$ is the map that sends
$a\ot b\ot c$ to $b\ot c\ot a$, and $\tau$ is the flip $a\ot b \mapsto
b\ot a$. We will show that $\a$ is a right weakly multiplicative
coaction, and that $\b$ is a left weakly multiplicative coaction.

\begin{lemma}\label{coactions}
In the above notation, the maps $\a$ and $\b$ are 
right coactions. 
\end{lemma}

\begin{proof} 
We do the computations for $\a$. The computations for $\b$ are similar.
First, to show that $\a$ is a coaction, we have to show that
$(\id\ot\bigd)\circ\a = (\a\ot\id)\circ\a$ and that $(\id\ot\aug)\circ\a
= \id_A$. The commutativity assumption on $\l$ and $\r$ make it possible
to use the Sweedler notation. Note that in this notation, the formula
for $\a$ is $\a(a) = \suma a_0\ot a_1S(a_{-1})$. Thus, 
$(\a\ot\id)\circ\a(a) = \suma a_0\ot a_1S(a_{-1})\ot a_2S(a_{-2})$. Now,
$(\id\ot\bigd)\circ\a(a) = (\id\ot\bigd)(\suma a_0\ot a_1S(a_{-1})) =
\suma a_0\ot \bigd(a_1S(a_{-1})) = \suma a_0\ot
\left\{\bigd(a_1)\cdot\bigd(S(a_{-1}))\right\} = \suma a_0\ot
\left\{(a_1\ot a_2)\cdot(S(a_{-1})\ot S(a_{-2}))\right\} =
\suma a_0\ot a_1S(a_{-1})\ot a_2S(a_{-2})$, as required. 

Next, $(\id\ot\aug)\circ\a(a) = (\id\ot\aug)(\suma a_0\ot a_1S(a_{-1}))
= \suma a_0\ot \aug(a_1)\aug(S(a_{-1})) = \suma a_0\ot
\aug(a_1)\aug(a_{-1}) = \suma \aug(a_{-1})a_0\ot\aug(a_{1})
=  \suma a_0\ot\aug(a_{1}) = \suma a_0\aug(a_{1}) = a$, as required.
Thus, $\a$ is a coaction. 
\end{proof}

\begin{proposition}\label{weakmult}     
In the above notation, the coaction $\a$ is right weakly multiplicative
and the coaction $\b$ is left  weakly multiplicative. 
In particular, the sets of coinvariants $A^{\coa}$ and
$A^{\cob}$ are subalgebras.
\end{proposition}

\begin{proof} The proofs are similar to \cite[Lemma 2.2]{dl} and
\cite[Proposition 1.1]{glr}. 
We give the proof for $\a$, the proof for $\b$ is similar. 
Let $a, b \in A$ and suppose that $\a(b) = b\ot 1$. Then, 
$\a(ab) = \sumab (ab)_0\ot (ab)_1S((ab)_{-1}) = \suma\sumb 
a_0b_0\ot a_1b_1 
S(b_{-1})S(a_{-1}) = \suma
(a_0\ot a_1)\cdot(\sumb b_0\ot b_1S(b_{-1}))\cdot(1\ot
S(a_{-1})) = 
\suma (a_0\ot a_1)\cdot(b\ot 1)\cdot(1\ot S(a_{-1})) = \suma (a_0\ot
a_1S(a_{-1}))\cdot(b\ot 1) = \a(a)\a(b)$; so that $\a$ is right weakly
multiplicative. Now, apply Lemma~\ref{subalg}. 
\end{proof}

In \cite{dl, gl, glr} coactions that are not algebra homorphisms are
constructed in this way, and their coinvariants are calculated.  
In \cite{dl}, for example, $A=\ch=\oglqn$, and $\alpha$, $\beta$ 
are two possible quantum deformations of the conjugation coaction. 
Explicit generators and the structure of the two subalgebras of 
coinvariants were determined for the case when $q$ is not a root 
of unity. The commutativity (and a vector space isomorphism) 
of these coinvariant subalgebras was explained recently in 
\cite{cw}; it was shown to 
follow from the fact that $\oglqn$ is coquasitriangular, 
and an interpretation of coinvariants via cocommutativity conditions
(see Section~\ref{frt}). 
However, in \cite{dl} (and in \cite{cw}) 
the connection between these two possible quantum deformations $\alpha,\beta$ 
of the conjugation coaction is not considered. 
This is done in the next section of this
paper for coquasitriangular Hopf algebras.

\section{The coquasitriangular case}\label{coquasi}
\marginpar{coquasi}

Recall, from \cite[10.1.1]{ks}, that a bialgebra $\ch$ is {\em
coquasitriangular} if 
there exists a linear form 
$\br:
\ch\ot\ch\goesto k$ such that the following conditions hold: \newline 
(i)~~$\br$ is convolution invertible with inverse $\brbar:
\ch\ot\ch\goesto k$; that is, $\br\star\brbar = \brbar\star\br = \aug\ot
\aug$ on $\ch\ot\ch$, \newline
(ii)~~$m^{{\rm op}} = \br\star m\star \brbar$ on $\ch\ot\ch$, \newline
(iii)~~$\br\circ(m\ot\id) = \br_{13}\star\br_{23}$ and $\br\circ(\id\ot m) =
\br_{13}\star\br_{12}$ on  $\ch\ot\ch\ot\ch$, \newline 
where $m^{{\rm op}}(a\ot b) = ba, m(a\ot b) = ab, \br_{12}(a\ot b\ot c)
= \br(a\ot b)\aug(c), \br_{23}(a\ot b\ot c) = \aug(a)(b\ot c)$ and
$\br_{13} (a\ot b\ot c) = \aug(b)(a\ot c)$, for $a,b,c \in \ch$.  
The map $\br$ is known as a {\em universal
$r$-form} of $\ch$.

In Sweedler notation, condition (i) is 
\[
\sum\, \br(\aone, \bone)\brbar(\atwo, \btwo) =
\sum\, \brbar(\aone, \bone)\br(\atwo, \btwo) = \aug(a)\aug(b); 
\]
and conditions (ii) and (iii) are 
\begin{align*}
ba &= \sum\, \br(\aone, \bone)\atwo\btwo\brbar(\athree, \bthree)\\
\br(ab,c) &= \sum\, \br(a, \cone)\br(b,\ctwo)\\
\br(a, bc) &= \sum\, \br(\aone, c)\br(\atwo, b)
\end{align*}
for $a, b, c\in \ch$. Here, we are thinking of $\br$ as a bilinear form
$\ch\times\ch\goesto k$, and $\br(a,b) = \br(a\ot b)$. 

\begin{proposition}\label{isocoact}
Let $\ch$ be a coquasitriangular Hopf algebra. Suppose that $A$ is an
algebra and that $\lambda: A \goesto \ch\ot A$ and $\rho:A\goesto
A\ot\ch$ are algebra homomorphisms that are left and right coactions,
respectively, and that $\lambda$ and $\rho$ commute. Suppose that $\a$
and $\b$ are defined as in the previous section. Then $\a$ and $\b$ are
isomorphic coactions. 
\end{proposition}

\begin{proof}
Let $\br$ be a universal $r$-form for $\ch$, with convolution inverse
$\brbar$.  Define $\psi:A\goesto 
A$ by $\psi(v) = \sumv\br(S(v_{-1}), v_1)v_0$. 

First we show that $\psi$ is a vector space isomorphism, by giving an
explicit inverse $\phi$. Define $\phi:A\goesto
A$ by $\phi(v) = \sumv\brbar(S(v_{-1}), v_1)v_0$. Then, 
\begin{align*}
\phi(\psi(v)) &= \phi(\sumv \br(S(v_{-1}), v_1)v_0)\\
	&= \sumv \br(S(v_{-2}), v_2)\brbar(S(v_{-1}), v_1)v_0\\
	&= \sumv \aug(S(v_{-1}))v_0\aug(v_1)\\
	&= v, 
\end{align*}
and, similarly, one checks that $\psi(\phi(v)) =v$; so that $\psi$ and
$\phi$ are inverse linear transformations. 

Next, we check that $\psi$ intertwines with $\a$ and $\b$; that is,
$\a\circ\psi = (\psi\ot\id)\circ\b$. An equivalent condition to 
(ii) above 
is the requirement that $\sum\,
\br(\aone, \bone)\atwo\btwo = \sum\,\br(\atwo, \btwo)\bone\aone$, for
$a, b \in \ch$, see \cite[p332, equation (5)]{ks}. Applying this with $a
= S(v_{-1})$ and $b = v$, and noting that 
$\Delta(S(v_{-1})) = \sum S(v_{-1}) \ot S(v_{-2})$, we obtain 
$\sum\, \br(S(v_{-1}), v_1)S(v_{-2})v_2 = \sum\, \br
(S(v_{-2}),v_2)v_1S(v_{-1})$. Thus,
\begin{align*}
\a(\psi(v)) &= \a(\sumv \br(S(v_{-1}), v_1)v_0)\\
  	&= \sumv \br(S(v_{-2}), v_2)v_0 \ot v_1S(v_{-1})\\
	&= \sumv v_0\ot \br(S(v_{-2}), v_2)v_1S(v_{-1})\\
	&= \sumv v_0\ot \br(S(v_{-1}), v_1)S(v_{-2})v_2\\
	&= \sumv \br(S(v_{-1}), v_1)v_0\ot S(v_{-2})v_2\\
	&= (\psi\ot\id)(\b(v))
\end{align*}
for $v\in A$; that is, $\a\circ\psi = (\psi\ot\id)\circ\b$, as required.
\end{proof} 

	It is now easy to check that $\b\circ\phi = (\phi\ot\id)\circ\a$
and that $v\in A$ is a $\b$-coinvariant if and only if $\psi(v)$ is an
$\a$-coinvariant. 

	In Section~\ref{frt}, 
we will see in examples that there is often a
natural grading present on the algebra $A$, and that this induces a
grading on the subalgebras of coinvariants for $\a$ and $\b$. The maps
$\psi$ and $\phi$ respect these gradings, and so the subalgebras of
coinvariants automatically have the same Hilbert series. However, the
maps $\psi$ and $\phi$ are not algebra homomorphisms in general 
(nor are their restrictions to the coinvariant subalgebras).

\section{Coinvariant subalgebras are finitely generated}\label{cofingen}
\marginpar{cofingen}

Next, we identify conditions that will guarantee that the coinvariants
of a weakly multiplicative coaction form a finitely generated
algebra.

Recall that a Hopf algebra $\ch$ is {\em cosemisimple} if any
$\ch$-comodule decomposes as the direct sum of irreducible subcomodules,
\cite[11.2.1]{ks}.  An equivalent condition for $\ch$ to be cosemisimple
is given in terms of Haar functionals.

\begin{definition} A linear functional $h$ on $\ch$ is {\em left
invariant} if $(\id \ot h)\Delta(a) = h(a)1$, for all $a\in \ch$.
Similarly, $\ch$ is {\em right invariant} if $(h \ot \id)\Delta(a) =
h(a)1$, for all $a\in \ch$.
\end{definition}

A Hopf algebra $\ch$ is {\em cosemisimple} if and only if there exists a
unique left and right invariant linear functional $h$ on $\ch$ such that
$h(1) =1$, \cite[11.2.1, Theorem 13]{ks}.  Such a functional is then
called the {\em Haar functional} of $\ch$. In Sweedler notation, the
Haar functional satisfies $h(a)1 = \sum_{(a)}\, a_1h(a_2) = \sum_{(a)}\,
h(a_1)a_2$ for all $a\in \ch$.

Our aim is to show that for suitable 
left or right weakly multiplicative right
coactions $\b$ of a cosemisimple Hopf algebra $\ch$ the coinvariants of
$\b$ form a finitely generated subalgebra. 

\begin{lemma}\label{haar}
Let $\ch$ be a cosemisimple Hopf algebra, with Haar
functional $h$, and let $V$ be a $k$-vector space such that 
$\b:V\goesto V\ot \ch$ is a right coaction. 
Set $\pi:= (\id \ot h)\circ\b:V\goesto V$. Then (i) $\pi(V) \sse V^{\cob}$,
(ii) $\pi(v) = v$ for all $v\in V^{\cob}$ and (iii) $\pi^2 = \pi$. 
\end{lemma} 

\begin{proof} 
(i) Let $v\in V$. Note that $\pi(v) = (\id \ot h)\circ\b(v) = \sum_{(v)}\,
v_0h(v_1)$. Hence, 
\begin{align*}
\b(\pi(v)) &= \b(\sum_{(v)}\, v_0h(v_1))=\sum_{(v)}\, v_0 \ot
v_1h(v_2)\\
	&=  \sum_{(v)}\, v_0 \ot h(v_1)1 
=  (\sum_{(v)}\, v_0h(v_1))\ot 1 \\
	&= \pi(v) \ot 1, 
\end{align*}
as required. (Note that the second equality is justified by the fact
that $\b\circ(\id \ot h)\circ \b = (\id\ot\id\ot
h)\circ(\b\ot\id)\circ\b = (\id\ot\id\ot h)\circ(\id\ot\Delta)\circ\b$.)

\noindent (ii) Let $v \in V^{\cob}$. Then 
$\pi(v)  = (\id \ot h)\circ\b(v) = (\id \ot h)\{v\ot 1\} = vh(1) = v$.

\noindent (iii) Let $v\in V$. Then $\pi(v) \in V^{\cob}$ by
(i); so the result follows by (ii).
\end{proof}

Thus, $\pi:= (\id \ot h)\circ\b$ is a projection from $V$ onto $V^{\cob}$,
cf. \cite[11.2.2, Corollary 19]{ks} and the comment following that
result.

\begin{lemma}\label{proj} 
Let $A$ be a $k$-algebra and $\ch$ be a cosemisimple Hopf algebra over
$k$.  Suppose that $\b:A \goesto A\ot \ch$ is a left weakly
multiplicative right coaction with $R:= A^{\cob}$.  Set $\pi:= (\id \ot
h)\circ\b$.  Then $\pi$ is a left $R$-module epimorphism from $A$ to
$R$; that is, $\pi(A) = R$ and 
$\pi(ba) = b\pi(a)$ for all $a\in A$ and $b\in R$. 
\end{lemma}

\begin{proof} 
Parts (i) and  (ii) of the previous lemma show 
that $\pi(A) = R$. Let $b\in
R$ and $a\in A$. Then 
\begin{align*}
\pi(ba) &= (\id \ot h)\{\b(ba)\} = (\id \ot h)\{\b(b)\b(a)\}\\
	&= (\id \ot h)\{(b\ot 1)(\sum_{(a)}\, a_0\ot a_1)\} = (\id \ot
h)\{(\sum_{(a)}\, ba_0 \ot a_1\}\\
	&= \sum_{(a)}\, ba_0h(a_1) = b\{\sum_{(a)}\, a_0h(a_1)\}\\
	&= b\pi(a), 
\end{align*}
so that $\pi$ is an $R$-module homomorphism. 
\end{proof}

The following proposition is well-known; we give a proof for
completeness. 

\begin{proposition} \label{well-known}
Suppose that $R = R_0 \oplus R_1 \oplus \dots$ is an
$\mn$-graded subring of a right noetherian ring $A$. Suppose that there
exist an epimorphism $\pi: A \goesto R$ of left $R$-modules, such that
$\pi(A) = R$ and  $\pi(R_0) = R_0$. Then there exist homogeneous
elements $r_1, \dots, r_s \in R$, of positive degree, such that $R =
\sum_{u\in M}\, uR_0 $, where $M$ is the multiplicative semigroup
generated by $1, r_1, \dots, r_s$. In particular, $R = R_0\langle r_1,
\dots, r_s \rangle$ as a ring. 
\end{proposition} 

\begin{proof} 
Set $R^+:=R_1 \oplus R_2 \oplus \dots$. Then $R^+A$ is a right ideal of
the right noetherian ring $A$; so there are homogeneous elements $r_i
\in R^+$, for $i=1, \dots, s$, with $R^+A = r_1A + \dots + r_sA$.
Observe that $R \sse R_0 + R^+A$. We
claim that $R = \sum_{u\in M}\, uR_0 $, where $M$ is the multiplicative
semigroup generated by $1, r_1, \dots, r_s$. Note that $R_0 \sse
\sum_{u\in M}\, uR_0 $. Assume that $n>0$ and that $R_m \sse \sum_{u\in
M}\, uR_0 $ for each $0\leq m < n$. Now 
\[
R_n \sse R = \pi(R) \sse \pi(R_0 +R^+A) = R_0 + \pi(\sum_{i=1}^{s}\,
r_iA) \sse R_0 + \sum_{i=1}^{s}\, r_i\pi(A) = R_0 + \sum_{i=1}^{s}\, r_i R.
\]
Let $r \in R_n$. Then there exist $r_0 \in R_0$ and $x_i \in R$ with $r=
r_0 + \sum_{i=1}^{s}\, r_ix_i$. Without loss of generality, we may
assume that the $x_i$ are homogeneous and that each $r_ix_i \in R_n$.
Thus, for each $i$, we have $x_i \in R_m$ for some $m<n$, since $r_i \in
A^+$. By the inductive hypothesis, 
each $x_i \in \sum_{u\in M}\, uR_0 $ and the
result follows. 
\end{proof} 

\begin{theorem}\label{hilbert}
Let $A = A_0 \oplus A_1 \oplus \dots $ be a right noetherian $\mn$-graded
$k$-algebra with $A_0 = k$ and let $\ch$ be a cosemisimple Hopf algebra
over $k$.  Suppose that $\b:A\goesto A\ot\ch$ is a left weakly multiplicative
right coaction such that $\b(A_i) \sse A_i\ot\ch$ for each $i$.  Then
$A^{\cob}$ is a subalgebra of $A$ that is finitely generated as a
$k$-algebra. 
\end{theorem}

\begin{proof}
Note that $R:= A^{\cob}$ is a subalgebra of $A$, by Lemma~\ref{subalg}. 
Set $R_i:= R\cap A_i$ and note that $R = R_0 \oplus R_1 \oplus \dots$, since
$\b(A_i)\sse A_i\ot\ch$.  Set $\pi := (\id\ot h)\circ\b$.  Then $\pi$ is
a left $R$-module epimorphism, by Lemma~\ref{proj}.  The result now
follows from Proposition~\ref{well-known}.  
\end{proof}

\begin{remark}\label{left-noetherian} 
{\rm The obvious modification of 
Theorem~\ref{hilbert} holds for a left noetherian algebra 
endowed with a right weakly multiplicative right coaction.} 
\end{remark}  

\section{Coinvariants of quantum groups in FRT-bialgebras
}\label{frt}
\marginpar{frt}

Examples of coactions $\alpha,\beta$ studied in Section~\ref{coinvs} 
are naturally associated to quantized algebras of functions on 
the simple Lie groups $SL(N)$, $O(N)$, $Sp(N)$. 
Let us briefly recall their construction due to 
Faddeev-Reshetikhin-Takhtajan \cite{rtf}; 
see \cite[Chapter 9]{ks} for a detailed presentation. 

Take the free associative algebra 
$\mc\langle x_{ij}\rangle$ in $N^2$ generators $x_{ij}$, $i,j=1,\ldots,N$. 
Fix an $N^2\times N^2$ matrix $R$ with complex entries 
$R_{ij}^{mn}$, $i,j,m,n=1,\ldots,N$. 
Let $\car$ be the quotient algebra of 
$\mc\langle x_{ij}\rangle$ modulo the two-sided ideal generated by 
\begin{equation}\label{eq:frt} 
\sum_{k,l}R^{ji}_{kl}x_{km}x_{ln}
=\sum_{k,l}x_{ik}x_{jl}R^{lk}_{mn}, 
\ \ i,j,m,n=1,\ldots,N. 
\end{equation} 
Since this ideal is generated by homogeneous elements 
(with respect to the standard grading of the free algebra), 
the algebra $\car$ is graded, the generators $x_{ij}\in\car$ 
having degree $1$. 
There is a unique bialgebra structure on $\car$ such that 
\begin{equation} 
\Delta(x_{ij})=\sum_kx_{ik}\otimes x_{kj}, \ \ 
\varepsilon(x_{ij})=\delta_{ij}, 
\end{equation} 
where $\delta_{ii}=1$, and $\delta_{ij}=0$ for $i\neq j$. 
Note that the comultiplication and the counit respect the grading.  
The bialgebra $\car$ is called the 
{\it coordinate algebra of the quantum matrix space} 
associated with $R$, or briefly, the {\it FRT-bialgebra}. 
It is shown in \cite{bg} that if 
$R$ is lower (or upper) triangular, then $\car$ is a noetherian 
algebra. 

For the rest of this section we assume that $q\in\mc^*$ is not a root of 
unity, and that 
$R$ is the R-matrix of the vector representation of one of the 
Drinfeld-Jimbo algebras 
$U_q(\slN)$, 
$U_q(\soN)$ with $N=2n$, 
$U_{q^{1/2}}(\soN)$ with $N=2n+1$, 
$U_q(\spN)$ with $N=2n$. 
See \cite[9.2, 9.3]{ks} for the explicit formulae for $R$; 
for the case of $\slN$, the relations \eqref{eq:frt} 
will be given in Section~\ref{gl-sl}, \eqref{definition-of-quantum-matrices}. 
For each of these cases $\ca:=\car$ contains a distinguished 
central group-like element $Q$ 
(for the case of $\slN$, $Q$ is the quantum determinant, 
and $Q$ is a quadratic element in the other cases, given 
explicitly in \cite[9.3.1]{ks}), 
such that the quotient bialgebra 
$\ch:=\car/\langle Q-1\rangle$ has a unique Hopf algebra structure. 
The Hopf algebra $\ch$ is denoted by 
$\oslqn$, $\ooqn$, $\ospqn$ in the respective cases, 
and is called the 
{\it coordinate algebra of the quantum 
special linear group, orthogonal group}, and 
{\it symplectic group}, respectively. 
Furthermore, $\ooqn$ contains a group-like element $D$ 
(similar to the quantum determinant). 
The quotient Hopf algebra $\ooqn/\langle D-1\rangle$ is denoted by 
$\osoqn$, and is called the 
{\it coordinate algebra of the quantum 
special orthogonal group} (cf. \cite{h}). 

From now on $\ch$ stands for any of 
$\oslqn$, $\ooqn$, $\osoqn$, $\ospqn$, 
and $\pi:\ca\to\ch$ denotes the natural homomorphism from 
the corresponding FRT-bialgebra. 
If $q$ is transcendental over $\mq$, then $\ch$ is cosemisimple 
by \cite{h}. The cosemisimplicity of $\oslqn$ is known under the 
weaker assumption that $q$ is not a root of unity, see \cite{t} 
(or \cite{nym} combined with our Section~\ref{gl-sl}). 
It is well known that the bialgebras $\ca$, 
$\oslqn$, $\ooqn$, $\ospqn$ are coquasitriangular, see \cite[10.1.2]{ks}. 

Set 
\[\lambda:=(\pi\otimes\id)\circ\Delta:\ca\to\ch\otimes\ca,\]
\[\rho:=(\id\otimes\pi)\circ\Delta:\ca\to\ca\otimes\ch.\] 
Clearly $\lambda$ is a left coaction, $\rho$ is a right coaction, 
and both $\lambda$ and $\rho$ are algebra homomorphisms. 
Moreover, $\lambda$ and $\rho$ commute 
by the coassociativity of $\Delta$. 
So we are in the setup of 
Section~\ref{coinvs}, and the formulae \eqref{eq:alpha}, \eqref{eq:beta} 
of Section~\ref{coinvs} 
define the right coactions $\alpha$ and $\beta$ of $\ch$ on $\ca$. 
First we observe that the 
coinvariants can be interpreted independently from the 
coactions $\alpha$, $\beta$ as follows. 

\begin{proposition} \label{pi-cocomm}
For $f\in \ca$, with $\Delta(f)=\sum f_1\otimes f_2$ we have 
\begin{itemize}
\item[(i)] $f$ is an $\alpha$-coinvariant if and only if 
$\sum f_1\otimes\pi(f_2)=\sum f_2\otimes \pi(f_1)$. 
\item[(ii)] $f$ is a $\beta$-coinvariant if and only if 
$\sum f_1\otimes f_2=\sum f_2\otimes S^2(\pi(f_1))$. 
\end{itemize}
\end{proposition} 

\begin{proof} 
This is a straightforward modification of the proof of 
\cite[Theorem 2.1]{dl} and \cite[Lemma 1.1]{cz}. 
\end{proof} 

As an application of our results, properties of the subalgebras of 
$\alpha$-coinvariants and $\beta$-coinvariants 
in $\ca$ are summarized in the following theorem. 

\begin{theorem}\label{main} 
Let $\ch$ be any of $\oslqn$, $\ooqn$, $\ospqn$, 
and let $\ca$ be the corresponding FRT-bialgebra. 
Then 
\begin{itemize} 
\item[(i)] the sets of coinvariants $\ca^{\coa}$ and $\ca^{\cob}$ 
are graded subalgebras of $\ca$; 
\item[(ii)] the algebras $\ca^{\coa}$ and $\ca^{\cob}$ are commutative. 
\end{itemize} 
Assume in addition that $q$ is transcendental over $\mq$ 
(or $\ch=\oslqn$ and $q$ is not a root of unity). Then 
\begin{itemize} 
\item[(iii)] $\ca^{\coa}$ and $\ca^{\cob}$ are finitely generated 
subalgebras of $\ca$;   
\item[(iv)] the graded algebras $\ca^{\coa}$ and $\ca^{\cob}$ 
have the same Hilbert series. 
\end{itemize} 
\end{theorem} 

\begin{proof} (i) The coinvariants form a subalgebra by 
Proposition~\ref{weakmult}. This is a graded subalgebra, since the 
comultiplication on $\ca$ is homogeneous. 

(ii) Both $\ca$ and $\ch$ are coquasitriangular bialgebras. 
Moreover, by \cite[10.1.2, Theorem 9]{ks} we know that there is a universal 
r-form $\br$ on $\ca$ which induces a universal r-form on $\ch$ 
denoted by the same symbol $\br$. Thus we have 
\begin{equation}\label{eq:0}
\br(x,y)=\br(\pi(x),\pi(y)) 
\mbox{ for }x,y\in\ca.
\end{equation}  
Assume that $a,b$ are $\gamma$-coinvariants, where $\gamma$ is either 
$\alpha$ or $\beta$. 
Let $T$ denote the identity operator on $\ch$ when $\gamma=\alpha$, 
and the automorphism $S^2$ of $\ch$ when $\gamma=\beta$. 
Note that 
\begin{equation}\label{eq:2} 
\br(T(x),T(y))=\br(x,y)
 \mbox{ for any }x,y\in\ch
\end{equation} 
by \cite[10.1.1, Proposition 2.(v)]{ks}. 
Then we have 
\[\sum a_1\otimes a_2=\sum a_2\otimes T(\pi(a_1)) \ 
\mbox{ and }\ 
\sum b_1\otimes b_2=\sum b_2\otimes T(\pi(b_1))\] 
by Proposition~\ref{pi-cocomm}.
Applying $(\pi\otimes\id\otimes\id)\circ(\Delta\otimes\id)$ 
to these equalitites we get 
\begin{align}\label{eq:1}
\sum \pi(a_1)\otimes a_2\otimes\pi(a_3)
&=\sum \pi(a_2)\otimes a_3\otimes T(\pi(a_1))\\
\label{eq:3} \sum \pi(b_1)\otimes b_2\otimes\pi(b_3)
&=\sum \pi(b_2)\otimes b_3\otimes T(\pi(b_1)). 
\end{align} 
We have  
\[ba = \sum\, \br(\aone, \bone)\atwo\btwo\brbar(\athree, \bthree)\]
by condition (ii) in the defining properties of $\br$ 
(see Section~\ref{coquasi}). By 
\eqref{eq:0}, \eqref{eq:1}, \eqref{eq:3}, 
\eqref{eq:2}, and the defining property 
(i) of $\br,\brbar$ in Section~\ref{coquasi},  
the right hand side of the above equality equals 
\begin{align*} 
\sum\, \br(\pi(a_1),\pi(b_1))a_2b_2\brbar(\pi(a_3),\pi(b_3))
&= \sum\, \br(\pi(a_2),\pi(b_2))a_3b_3\brbar(T(\pi(a_1)),T(\pi(b_1))
\\
&=\sum\, \br( (a_2,b_2)a_3b_3 \brbar(a_1,b_1)   
\\
&=\sum\, \aug(a_1)\aug(b_1)a_2b_2=ab, 
\end{align*} 
thus $ab=ba$ for all $\gamma$-coinvariants $a,b$. 

(iii) The additional assumptions on $q$ 
ensure that $\ch$ is cosemisimple. As we noted above, $\ca$ is both left and 
right Noetherian by \cite{bg}, since the R-matrix is triangular. 
The homogeneous components of $\ca$ are subcomodules with respect to $\alpha$ 
and $\beta$. Hence we may apply Theorem~\ref{hilbert} to conclude 
that $\ca^{\coa}$ and $\ca^{\cob}$ are finitely generated graded subalgebras. 

(iv) $\alpha$ and $\beta$ are isomorphic coactions by 
Proposition~\ref{isocoact}; the map $\psi$ in its proof is homogeneous on 
$\ca$, hence the restrictions $\alpha_d$, $\beta_d$ of $\alpha$, $\beta$ 
to the homogeneous component $\ca_d$ of $\ca$ are also isomorphic coactions. 
The cosemisimplicity of $\ch$ implies that the dimension of the space of
$\alpha$-coinvariants in 
$\ca_d$ equals the multiplicity of the trivial corepresentation in 
$\alpha_d\cong\beta_d$, which equals the dimension of the 
space  of $\beta$-coinvariants in $\ca_d$. 
In other words, the graded algebras 
$\ca^{\coa}$ and $\ca^{\cob}$ have the same Hilbert series. 
\end{proof} 

\begin{remark}\label{remark-soqn}
{\rm The statements (i) and (iii) of Theorem~\ref{main} (and their proofs) 
are clearly valid for $\ch=\osoqn$ and the corresponding FRT-bialgebra.}    
\end{remark} 

\begin{remark}\label{remark-dlcw}
{\rm The above theorem was motivated by the main result 
of \cite{dl}. In that paper $\ch=\oglqn$, 
the coordinate algebra of the quantum 
general linear group, $\ca=\omqn$, the coordinate algebra of $N\times N$ 
quantum matrices, and $\pi$ is the natural embedding of $\omqn$ into 
$\oglqn$. The main result of \cite{dl} is that if $q$ is not a root of 
unity, then $\ca^{\coa}$ and $\ca^{\cob}$ are $N$-variable commutative 
polynomial algebras, with explicitly given generators. 
(This result implies the case of $\oslqn$ of Theorem~\ref{main}, 
as we shall show in Theorem~\ref{sl-coinv}.) Moreover, 
$\alpha$ and $\beta$ extend to a coaction of the Hopf algebra 
$\oglqn$ on itself, and the subalgebras of coinvariants can be explicitly 
described. It was observed later in \cite{cw} that the commutativity 
of these coinvariant subalgebras follows from the coquasitriangularity 
of $\oglqn$, and the interpretation of the coinvariants in terms of 
certain cocommutativity conditions. 
The proof of statement (ii) in Theorem~\ref{main} 
is a modification of the proof of \cite[Proposition 2.1]{cw}. 
} 
\end{remark}

\begin{remark}\label{remark-other-examples} 
{\rm Theorem~\ref{hilbert} can be applied to conclude the finite generation 
property of various other subalgebras of coinvariants arising in this context. 
For example, the algebra $\ca^{\col}$ of $\lambda$-coinvariants 
and the algebra $\ca^{\cor}$ of $\rho$-coinvariants are 
also finitely generated, provided that $\ch$ is cosemisimple. 
See also the problem studied in \cite{glr}. }
\end{remark} 

There are unique right coactions $\bar\alpha$, $\bar\beta$ 
of the Hopf algebra $\ch$ on itself such that $\pi$ intertwines 
between $\alpha$ and $\bar\alpha$, 
and $\pi$ intertwines between $\beta$ and $\bar\beta$. 
In Sweedler's notation we have 
\[\bar\beta(x)=\sum x_2\otimes S(x_1)x_2 \ 
\mbox{ and }\ 
\bar\alpha(x)=\sum x_2\otimes x_3S(x_1).\] 
Since $\bar\beta$ is the formal dual of the adjoint action
(see \cite[p.36]{m}), it is called 
the {\it adjoint coaction} of $\ch$; the version 
$\bar\alpha$ was introduced in \cite{dl}. 
Recall that $x\in\ch$ is {\it cocommutative} if 
$\sum x_1\otimes x_2=\sum x_2\otimes x_1$, 
Following \cite{cw}, we call $x\in\ch$ 
$S^2$-{\it cocommutative} if 
$\sum x_1\otimes x_2=\sum x_2\otimes S^2(x_1)$. 
The $S^2$-cocommutative elements coincide with the $\bar\beta$-coinvariants 
(see \cite[Lemma 1.1]{cz}, 
whereas the cocommutative elements coincide with the 
$\bar\alpha$-coinvariants (see \cite[Theorem 2.1]{dl}). 

\begin{corollary} \label{fingen-cocomm}
Let $\ch$ be any of $\oslqn$, $\ooqn$, $\osoqn$, $\ospqn$, 
and assume that $q$ is transcendental over $\mq$ (or 
$\ch=\oslqn$ and $q$ is not a root of unity). 
Then each of the cocommutative elements and the $S^2$-cocommutative 
elements in $\ch$ forms a finitely generated subalgebra. 
\end{corollary} 

\begin{proof} 
Since $\ch$ is cosemisimple by our assumption on $q$, the subspace 
of $\alpha$-coinvariants (respectively $\beta$-coinvariants) in $\ca$ 
is mapped onto the subspace of $\bar\alpha$-coinvariants 
(respectively $\bar\beta$-coinvariants) under 
the morphism $\pi$ of $\ch$-comodules. Moreover, 
$\pi$ is an algebra homomorphism. So the statement follows from 
Theorem~\ref{main} and the above remark. 
\end{proof} 
 
\begin{remark} {\rm 
It was proved in \cite[Theorem 2.2]{cw} that the 
coquasitriangularity of $\ch$ 
implies that the two subalgebras in Corollary~\ref{fingen-cocomm} 
are commutative. A vector space isomorphism between them was also 
established there. }
\end{remark}


\section{$\oglqn$ is a central extension of $\oslqn$}\label{gl-sl}
\marginpar{gl-sl} 

Over an algebraically closed base field, 
any element of the general linear group $GL(N)$ can be written as 
the product of a scalar matrix and an element from $SL(N)$. 
This and the centrality of scalar matrices imply a strong relationship  
between 
the representation  theories of $GL(N)$ and $SL(N)$. 
Our aim is to establish such a link between the corepresentation theories 
of $\oglqn$ and $\oslqn$. 
By a {\it corepresentation} we mean a right coaction of a Hopf algebra 
on a vector space. 

We work now over an arbitrary base field $k$, and $q$ is a non-zero 
element of $k$. 
Recall that the FRT-bialgebra of $\oslqn$ is $\omqn$, 
the {\em coordinate algebra of quantum $N\times N$ matrices}, which is the  
$k$-algebra generated by $N^2$ indeterminates $x_{ij}$, for 
$i,j=1,\dots,N$, subject to the following relations. 
\begin{equation} \label{definition-of-quantum-matrices}
\begin{array}{rcl}
x_{ij}x_{il} &=& qx_{il}x_{ij}, \\
x_{ij}x_{kj} &=& qx_{kj}x_{ij}, \\
x_{il}x_{kj} &=& x_{kj}x_{il},  \\
x_{ij}x_{kl} - x_{kl}x_{ij} &=& (q-q^{-1})x_{il}x_{kj},
\end{array}
\end{equation}
for $1 \le i < k \le N$ and $1 \le j < l \le N$. 
The algebra $\omqn$ is an iterated Ore extension, and so a noetherian domain. 
The {\em quantum determinant}, $\detq$, is the element 
\[
\detq := \sum_{\sigma \in S_{N}}(-q)^{l(\sigma)}
x_{1,\sigma(1)} \dots x_{N,\sigma(N)}.
\]
The element $\detq$ is central in $\omqn$ 
(see, for example, \cite[Theorem 4.6.1]{pw}), and by adjoining its 
inverse we get the 
{\em coordinate algebra of the quantum general linear group} 
\[
\oglqn := \omqn[\detq^{-1}].
\]
The algebra $\oglqn$ is a Hopf algebra, the {\em comultiplication} $\Delta$ 
being given by 
$\Delta(x_{ij})=\sum_{k=1}^Nx_{ik}\otimes x_{kj}$. 
The ideal generated by $(\detq-1)$ is a Hopf ideal of $\oglqn$.  
The {\em coordinate algebra of the quantum special linear group} 
$\oslqn$ is defined as the quotient Hopf algebra of $\oglqn$ modulo the  
ideal generated by $(\detq-1)$. 
Denote by $\pi$ the natural homomorphism 
$\oglqn\to\oslqn$, and set $y_{ij}:=\pi(x_{ij})$. 

The coordinate algebra of the multiplicative group of $k$ is 
the algebra of Laurent polynomials $k[z,z^{-1}]$.  
It is another Hopf algebra 
homomorphic image of $\oglqn$: the homomorphism 
$\pi_Z:\oglqn\to \ok$ is given by 
$\pi_Z(x_{ii})=z$ for $i=1,\ldots,N$ and 
$\pi_Z(x_{ij})=0$ if $i\neq j$. 
This can be paraphrased by saying that  
the multiplicative group of $k$ is 
a {\it quantum subgroup} of the quantum general linear group. 
Moreover, it is a {\em central quantum subgroup} in the following sense: 
\begin{equation}\label{eq:central-subgroup} 
\tau\circ(\pi_Z\otimes\id)\circ\Delta
=(\id\otimes\pi_Z)\circ\Delta
\end{equation}
where $\tau:\ok\otimes\oglqn\to\oglqn\otimes\ok$ is the flip 
$\tau(a\otimes b)=b\otimes a$. 

Consider $\eta:\ok\to\ozeta$, 
the natural homomorphism with $\eta(z)=\zeta$ and kernel 
$\langle z^N-1\rangle$. 
It follows from $\pi_Z(\detq)=z^N$ that 
$\eta\circ\pi_Z$ factors through $\pi$, so we have a commutative  
diagram 
\begin{equation}\label{eq:diagram} 
\begin{array}{ccc} 
\oglqn &\stackrel{\pi_Z}\longrightarrow &\ok\\
\downarrow{\pi} & &\downarrow{\eta}\\
\oslqn &\stackrel{\pi_{\zeta}}\longrightarrow &\ozeta
\end{array} 
\end{equation} 
of Hopf algebra homomorphisms. 
The map $\pi_{\zeta}$ is given explicitly by $\pi_{\zeta}(y_{ii})=\zeta$ 
for $i=1,\ldots,N$ and $\pi_{\zeta}(y_{ij})=0$ for $i\neq j$. 
Equation \eqref{eq:central-subgroup} implies 
\begin{equation}\label{eq:central-subgroup2} 
\tau\circ(\pi_{\zeta}\otimes\id)\circ\Delta
=(\id\otimes\pi_{\zeta})\circ\Delta .
\end{equation}

The relevance of the properties 
\eqref{eq:central-subgroup} and \eqref{eq:central-subgroup2} 
is shown by Proposition~\ref{homog-decomp} below. 
We state it in a general form. 
Let $\kappa:\ch\to\cc$ be a homomorphism of Hopf algebras such that 
\begin{equation}\label{eq:central}
\tau\circ(\id\otimes\kappa)\circ\Delta_{\ch}=
(\kappa\otimes\id)\circ\Delta_{\ch}. 
\end{equation}
Assume that $\cc$ is cosemisimple. Then $\cc$ decomposes as a direct sum 
$\bigoplus_{\varphi\in\Lambda}\cc(\varphi)$ 
of simple subcoalgebras, where 
$\Lambda$ is a complete list of irreducible corepresentations of $\cc$, 
and $\cc(\varphi)$ is the coefficient space of $\varphi$ (cf. \cite{g}). 
Any $\ch$-corepresentation $\psi:V\to V\otimes\ch$ 
`restricts' to a $\cc$-corepresentation 
$(\id\otimes\kappa)\circ\psi:V\to V\otimes\cc$. 
Therefore $V$ decomposes as 
$\bigoplus_{\varphi\in\Lambda}V_{\varphi}$, 
where 
\[V_{\varphi}=\{v\in V\mid 
(\id\otimes\kappa)\circ\psi(v)\in V\otimes\cc(\varphi)\}.\]
In other words, $V_{\varphi}$ is the sum of simple $\cc$-subcomodules of $V$ 
on which the $\cc$-corepresentation is isomorphic to $\varphi$. 
In particular, 
$\ch=\bigoplus_{\varphi\in\Lambda}\ch_{\varphi}$, 
where 
\[\ch_{\varphi}=\{f\in\ch\mid (\id\otimes\kappa)\circ\Delta(f)\in 
\ch\otimes\cc(\varphi)\}.\] 

\begin{proposition}\label{homog-decomp} 
Any summand $V_{\varphi}$ in $V=\bigoplus_{\varphi\in\Lambda}V_{\varphi}$ 
is an $\ch$-subcomodule, and 
$\psi(V_{\phi})\subseteq V_{\varphi}\otimes\ch_{\varphi}$. 
In particular, $\ch_{\varphi}$ is a subcoalgebra, and 
$V_{\varphi}$ is an $\ch_{\varphi}$-comodule. 
\end{proposition} 

\begin{proof} 
First we show that $\psi(V_{\varphi})\subseteq V\otimes\ch_{\varphi}$. 
Assume that $v\in V_{\varphi}$.  
Write 
$\psi(v)=\sum_i a_i\otimes b_i$ 
with $\{a_i\}$ linearly independent, 
and 
$\Delta(b_i)=\sum_jb_i^j\otimes c_i^j$. 
Then 
\begin{align*}
\sum_{i,j}a_i\otimes b_i^j\otimes\kappa(c_i^j)
&=(\id\otimes\id\otimes\kappa)\circ(\id\otimes\Delta)\circ\psi(v)
\\&=(\id\otimes\id\otimes\kappa)\circ(\psi\otimes\id)\circ\psi(v)
\\&=(\psi\otimes\id)\circ(\id\otimes\kappa)\circ\psi(v)
\end{align*}  
is contained in $V\otimes\ch\otimes\cc(\varphi)$. 
It follows that 
$\sum_jb_i^j\otimes \kappa(c_i^j)$ is contained in 
$\ch\otimes\cc(\varphi)$ for all $i$; 
that is, $b_i\in\ch_{\varphi}$ for all $i$. 
This means that $\psi(v)\in V\otimes\ch_{\varphi}$. 

So for $v\in V_{\varphi}$, 
\begin{equation}\label{eq:u}
u:=(\id\otimes\id\otimes\kappa)\circ(\id\otimes\Delta)\circ\psi(v)
\in V\otimes\ch_{\varphi}\otimes\cc(\varphi). 
\end{equation} 
Write $\psi(v)$ as a sum 
$\psi(v)=\sum_{\nu\in\Lambda}w_{\nu}$ 
with $w_{\nu}\in V_{\nu}\otimes\ch_{\varphi}$. 
By \eqref{eq:central} we have 
\begin{align*}
(\id\otimes\id\otimes\kappa)\circ(\id\otimes\Delta)\circ\psi
&=(\id\otimes\tau)\circ(\id\otimes\kappa\otimes\id)\circ
(\id\otimes\Delta)\circ\psi
\\&=(\id\otimes\tau)\circ
(\id\otimes\kappa\otimes\id)\circ(\psi\otimes\id)\circ\psi
\\&=(\id\otimes\tau)\circ
[((\id\otimes\kappa)\circ\psi)\otimes\id]\circ\psi. 
\end{align*} 
Thus $u=\sum_{\nu\in\Lambda}u_{\nu}$, where 
\begin{equation}\label{eq:unu}
u_{\nu}=(\id\otimes\tau)\circ
[((\id\otimes\kappa)\circ\psi)\otimes\id]\,(w_{\nu})\in
V_{\nu}\otimes\ch_{\varphi}\otimes\cc(\nu). 
\end{equation} 
By combining \eqref{eq:u} and \eqref{eq:unu} we conclude 
that 
$u_{\nu}=0$ if $\nu\neq\varphi$, hence $w_{\nu}=0$ if $\nu\neq\varphi$. 
That is, 
$\psi(v)=w_{\varphi}\in V_{\varphi}\otimes\ch_{\varphi}$, 
as we claimed. 

If we apply the result just proved for the corepresentation $\Delta$ of 
$\ch$, we get that $\ch_{\varphi}$ is a subcoalgebra. 
\end{proof} 

The Hopf algebra $\ok$ is the sum of one-dimensional subcoalgebras, 
hence it is cosemisimple, and 
all its irreducible corepresentations are $1$-dimensional. 
Up to isomorphism, the irreducible corepresentations are 
$k\to k\otimes\ok$, $1\mapsto 1\otimes z^d$, as $d$ ranges over the set of 
integers. 
Apply Proposition~\ref{homog-decomp} with 
$\ch=\oglqn$, $\cc=\ok$, and $\kappa=\pi_Z$. 
We get 
\[
\oglqn=\bigoplus_{d\in\mz}\oglqn_d,  
\] 
where 
\[\oglqn_d=\{f\in\oglqn\mid
(\id\otimes\pi_Z)\circ\Delta (f)=f\otimes z^d\}\] 
is a subcoalgebra of $\oglqn$. For an arbitrary corepresentation 
$\psi:V\to V\otimes\oglqn$ 
we have 
\begin{equation}\label{eq:homog-decomp}
V=\bigoplus_{d\in\mz}V_d, \mbox{ where }
V_d:=\{v\in V\mid (\id\otimes\pi_Z)\circ\psi(v)=v\otimes z^d\}.
\end{equation}  
By Proposition~\ref{homog-decomp}, $V_d$ is a subcomodule, moreover, 
it is an $\oglqn_d$-comodule.  

Similarly, $\ozeta$ is cosemisimple, its irreducible corepresentations are 
indexed by $\mz/(N)$, the set of residue classes modulo $N$: 
with $\bar d\in \mz/(N)$, the residue class of $d$ modulo $N$, 
we associate the one-dimensional corepresentation 
$1\mapsto 1\otimes \zeta^d$. 
Apply Proposition~\ref{homog-decomp} 
with $\ch=\oslqn$, $\cc=\ozeta$, and $\kappa=\pi_{\zeta}$. 
So 
\[
\oslqn=\bigoplus_{d\in\mz/(N)}\oslqn_{\bar d}, 
\]
where  
\[\oslqn_{\bar d}=\{f\in\oslqn\mid 
(\id\otimes \pi_{\zeta})\circ\Delta (f)=f\otimes \zeta^d\}\]
is a subcoalgebra of $\oslqn$. 
Any corepresentation $\psi:V\to V\otimes\oslqn$ decomposes as
\[
V=\bigoplus_{\bar d\in\mz/(N)}V_{\bar d}, 
\mbox{ where }
V_{\bar d}:=\{v\in V\mid (\id\otimes\pi_{\zeta})\circ\psi(v)
=v\otimes \zeta^d \}. 
\]
By Proposition~\ref{homog-decomp}, 
$V_{\bar d}$ is a subcomodule, moreover, it is an 
$\oslqn_{\bar d}$-comodule. 

Let us introduce some terminology, based on 
on the above discussion. 
We say that an $\oglqn$-corepresentation $\psi$ 
on $V$ is {\it homogeneous of degree} $d\in\mz$ if $V=V_d$;  
in this case $V$ is an $\oglqn_d$-comodule. 
Similarly, we say that an $\oslqn$-corepresentation $\psi$ on $V$
is {\it homogeneous of degree} $e\in\mz/(N)$ if $V=V_e$; 
in this case $V$ is an $\oslqn_e$-comodule. 
Note that an indecomposable corepresentation of $\oglqn$ 
or $\oslqn$ is necessarily homogeneous. 

Denote by $\pi_d$ the restriction of $\pi$ to $\oglqn_d$. 
The diagram \eqref{eq:diagram} implies that $\pi_d$ maps 
$\oglqn_d$ into $\oslqn_{\bar d}$. Our main aim is to show that $\pi_d$ is 
an isomorphism (of coalgebras) for all $d$. 
To see $\oglqn_d$ more explicitly, note that 
for the generators of $\oglqn$ we have 
$x_{ij}\in\oglqn_1$, $\detq^{-1}\in\oglqn_{-N}$, hence $\oglqn_d$ is 
spanned by the elements $w/\detq^{-r}$, where $w$ is a monomial in the 
$x_{ij}$ of degree $rN+d$. 
The space $\oslqn_{\bar d}$ is spanned by the monomials in the $y_{ij}$ 
whose degree is congruent to $d$ modulo $N$. 

The tensor product of Hopf algebras becomes a Hopf algebra with obvious 
operations. Consider 
\[\cb:=\oslqn\otimes\ok.\] 
We shall identify $\oglqn$ with a Hopf subalgebra of $\cb$. 
The technical advantage of $\cb$ is that the $\mz$-grading 
\[\cb=\bigoplus_{d\in\mz}\cb_d=\bigoplus_{d\in\mz}\oslqn\otimes z^d\] 
becomes visible. We have 
$\cb_d\cdot\cb_m=\cb_{d+m}$, so $\cb$ is a graded algebra. 
From $\Delta(z)=z\otimes z$ we get that 
$\Delta(\cb_d)\subseteq \cb_d\otimes \cb_d$, 
so $\cb_d$ is a subcoalgebra of $\cb$. 
Moreover, we have the obvious Hopf algebra surjections 
\[\id\otimes\varepsilon:\cb=\oslqn\otimes\ok\to \oslqn\] 
and 
\[\varepsilon\otimes\id:\cb=\oslqn\otimes\ok\to\ok,\] 
where $\varepsilon$ denotes the counit map of the appropriate Hopf algebra. 
The second surjection can be used to define the right coaction 
$(\id_{\cb}\otimes(\varepsilon\otimes\id))\circ\Delta_{\cb}:
\cb\to\cb\otimes\ok$, 
and the grading on $\cb$ can be reinterpreted in terms of this 
$\ok$-coaction as in \eqref{eq:homog-decomp}: 
\begin{equation}\label{eq:grading}
\cb_d=\{a\in\cb\mid (\id_{\cb}\otimes(\varepsilon\otimes\id))\circ
\Delta_{\cb}(a)=a\otimes z^d\}.
\end{equation}

\begin{proposition}\label{embed-glq} 
The map 
$\iota:=\pi\otimes\pi_Z:\oglqn\to\oslqn\otimes\ok$ 
is an embedding of the Hopf algebra 
$\oglqn$ into $\cb$. 
\end{proposition} 

\begin{proof} 
Both $\pi$ and $\pi_Z$ are homomorphisms of Hopf algebras, hence so is 
$\pi\otimes\pi_Z$. 

The only thing left to show is that $\iota$ is injective. 
Assume to the contrary that $f$ is a nonzero element of $\ker(\iota)$. 
Multiplying $f$ by an appropriate power of $\detq$ we may assume 
that $f\in\omqn$. 
Recall that 
$\cb$ is a $\mz$-graded algebra, and $\omqn$ is also $\mz$-graded,   
the generators $x_{ij}$ having degree $1$. 
Since $\iota(x_{ij})=y_{ij}\otimes z$, the restriction of 
$\iota$ to $\omqn$ is homogeneous. So 
all of the homogeneous components of $f$ 
are contained in $\ker(\iota)$. We may assume that $f$ itself 
is homogeneous of degree $d$. Then 
$0=\iota(f)=\pi(f)\otimes z^d$, implying that $\pi(f)=0$.  
Thus $f$ is a multiple of $(\detq-1)$. This is a contradiction, because no 
non-zero multiple of $(\det_q-1)$ in the domain $\omqn$ is homogeneous. 
\end{proof} 

\begin{remark}\label{levasseur-stafford}
{\rm Clearly $\iota(\oglqn)$ is the subalgebra of 
$\cb$ generated by 
$y_{ij}\otimes z$ ($1\leq i,j\leq N$) and $1\otimes z^{-N}$. 
It is easy to see that $\cb$ is a free $\oglqn$-module of rank $N$ 
generated by the central elements $1\otimes z^i$, $i=0,\ldots,N-1$, 
but we shall not use this fact. 
As $k$-algebras, $\oglqn$ and $\cb$ are isomorphic by \cite{ls}, 
where ring theoretic properties of $\oslqn$ are derived from this 
isomorphism. However, the algebra isomorphism $\oglqn\to\oslqn\otimes\ok$ 
given in \cite{ls} is not a morphism of coalgebras. 
To relate the corepresentation theories, we need the Hopf algebra 
morphism $\iota$. 
}
\end{remark} 

From now on we shall freely identify $\oglqn$ with 
$\iota(\oglqn)$. 
It is easy to check on the generators that $\pi$ is the restriction to 
$\oglqn$ of the map $\id\otimes\varepsilon$, and 
$\pi_Z$ is the restriction to $\oglqn$ of $\varepsilon\otimes\id$. 
It follows from \eqref{eq:grading} that 
\[\oglqn_d=\cb_d\cap\oglqn.\] 

\begin{proposition}\label{coalg-iso} 
The map $\pi_d$ is an isomorphism between the coalgebras 
$\oglqn_d$ and $\oslqn_{\bar d}$. 
\end{proposition} 

\begin{proof} By \eqref{eq:diagram} we know that $\pi_d$ is a surjective 
coalgebra homomorphism. So we only need to show the injectivity. 
The restriction of $\id\otimes\varepsilon$ to $\cb_d$ is clearly a vector 
space isomorphism $\cb_d\to\oslqn$. Since $\pi_d$ is the 
restriction of $\id\otimes\varepsilon$ to $\oglqn_d$, it is 
also injective. 
\end{proof} 

Let us record the following trivial observation. 

\begin{lemma}\label{subcoalg}
Let $D$ be a subcoalgebra of a coalgebra $C$, and 
let $\varphi:V\to V\otimes C$ be a corepresentation of $C$ such that 
$\im(\varphi)\subseteq V\otimes D$, so $V$ can be viewed as a $D$-comodule. 
Then a subspace $W$ of $V$ is a $C$-subcomodule if and only if $W$ is a 
$D$-subcomodule. 
In particular, $V$ is an irreducible (respectively indecomposable) 
$C$-comodule if and only if it is an 
irreducible (respectively indecomposable) $D$-comodule. 
\end{lemma}

\begin{corollary}\label{restrict-lift} 
(i) An indecomposable $\oglqn$-corepresentation $\psi$ of degree $d$ 
restricts to an indecomposable $\oslqn$-corepresentation 
$(\id\otimes\pi)\circ\psi$ of degree $\bar d$. 

(ii) Let $\varphi$ be an indecomposable $\oslqn$-corepresentation 
of degree $e\in\mz/(N)$. For $d\in \mz$ with $\bar d=e$, define 
$\varphi^{(d)}:=(\id\otimes\pi_d^{-1})\circ\varphi$. 
Then $\varphi^{(d)}$ is an indecomposable $\oglqn$-corepresentation of degree 
$d$, such that its restriction $(\id\otimes\pi)\circ\varphi^{(d)}$ is 
$\varphi$. 
\end{corollary} 

\begin{proof} This is an immediate consequence of Proposition~\ref{coalg-iso} 
and Lemma~\ref{subcoalg}. 
\end{proof} 

In particular, by Corollary~\ref{restrict-lift} 
any $\oslqn$-corepresentation $\varphi:V\to V\otimes\oslqn$ 
can be `lifted' to an $\oglqn$-corepresentation 
$\psi:V\to V\otimes \oglqn$, such that the restriction of 
$\psi$ to $\oslqn$ is $\varphi$ 
(note that $\psi$ is not uniquely determined). 

\begin{proposition}\label{irred-sl-gl}
Let $\psi:V\to V\otimes\oglqn$ be a corepresentation of $\oglqn$, 
and consider its restriction 
$\bar\psi:=(\id\otimes\pi)\circ\psi$ to $\oslqn$. 
\begin{itemize} 
\item[(i)]  
$\psi$ is irreducible if and only if 
$\bar\psi$ is irreducible. 
\item[(ii)] 
$\psi$ is cosemisimple if and only if 
$\bar\psi$ is cosemisimple. 
\end{itemize}
\end{proposition} 

\begin{proof} 
(i) It is trivial that the irreducibility of 
$\bar\psi$ implies the irreducibility of $\psi$. 
For the converse, observe that if $\psi:V\to V\otimes \oglqn$ is 
irreducible,  
then it is indecomposable, hence homogeneous of some degree $d$. 
In particular, it is an $\oglqn_d$-comodule, and the result follows 
using the coalgebra isomorphism $\pi_d$ of Proposition~\ref{coalg-iso},  
and Lemma~\ref{subcoalg}. 

(ii) Clearly $\psi$ is cosemisimple (that is, $V$ is the sum of 
simple subcomodules) if and only if its homogeneous components 
are cosemisimple. Thus we may assume that $\psi$ itself is homogeneous. 
Then the result immediately follows from the coalgebra isomorphism 
in Proposition~\ref{coalg-iso} and Lemma~\ref{subcoalg}, as above.  
\end{proof} 

Proposition~\ref{irred-sl-gl} (ii) gives an alternative proof 
of the fact that $\oglqn$ is cosemisimple if and only if 
$\oslqn$ is cosemisimple. Actually, both of them are known to be 
cosemisimple if and only if $q$ is not a root of unity, 
see \cite{nym}, \cite{h}, \cite{t}. 

Next we draw a consequence on $\oslqn$-coinvariants, 
extending to the quantum case 
the well known relation between relative $GL(N)$-invariants 
and absolute $SL(N)$-invariants. 

\begin{proposition}\label{rel-inv} 
Let $\psi$ be an $\oglqn$-corepresentation on $V$, and let 
$\bar\psi:=(\id\otimes\pi)\circ\psi$ be its restriction to $\oslqn$. 
Then the subspace of $\oslqn$-coinvariants in $V$ equals 
\[\bigoplus_{j\in\mz}\{v\in V\mid \psi(v)=v\otimes\detq^j\}.\]
In particular, if $\psi$ is homogeneous of degree zero, 
then $v\in V$ is an $\oslqn$-coinvariant if and only if 
$v$ is an $\oglqn$-coinvariant. 
\end{proposition} 

\begin{proof} Consider the decomposition $V=\bigoplus_{d\in\mz} V_d$ 
into the direct sum of homogeneous $\oglqn$-comodules. 
Obviously each summand $V_d$ is an $\oslqn$-subcomodule (with respect to 
$\bar\psi$). Therefore $v\in V$ is an $\oslqn$-coinvariant if and only if 
all the homogeneous components of $v$ are $\oslqn$-coinvariants. 
Thus we may restrict to the case when $V=V_d$ for some $d$. 
If there is a non-zero 
$\oslqn$-coinvariant in $V_d$, then $d=jN$ for an integer $j$, 
because $1\in \oslqn_{\bar 0}$. 
Assume that this is the case. 
We have $\pi_d(\detq^j)=1$, and $\pi_d$ is a coalgebra isomorphism 
by Proposition~\ref{coalg-iso}. 
So the coefficient space (cf. \cite{g}) of the $\oglqn$-comodule generated by 
$v$ is spanned by $\detq^j$. 
Moreover, we have $\psi(v)=v\otimes \detq^j$. 
\end{proof} 

For the rest of the paper we assume that $k=\mc$, the field of complex 
numbers. 
As an application of Proposition~\ref{rel-inv}, 
we deduce from the results of \cite{dl}  an explicit 
description of the subsets of $\oslqn$-coinvariants in 
$\omqn$ and $\oslqn$ 
with respect to the coactions $\alpha,\beta$ defined in 
Section~\ref{frt}. 

Fix an integer $t$ with $1\leq t \leq N$.  Let $I$ and $J$ be subsets of
$\{1,\dots, N\}$ with $|I|= |J| = t$.  The subalgebra of $\omqn$
generated by $x_{ij}$ with $i\in I$ and $j\in J$ can be regarded as an
algebra of $t\times t$ quantum matrices, and so we can calculate its
quantum determinant - this is a {\em $t\times t$ quantum minor} and we
denote it by $[I|J]$. 
The quantum minor $[I|I]$ is said to be a {\em principal quantum minor}.  
We denote the sum of all the principal
quantum minors of a given size $i$ by $\si$.  Note that $\sigma_1 =
x_{11} + \dots + x_{NN}$ and that $\sigma_N = \detq$.  
Consider also 
the {\em weighted sums of principal minors} 
$\ti:= \sum_{I}\, q^{-2w(I)}[I|I]$, 
$i=1,\ldots,N$ 
(here 
$w(I)$ denotes the sum of the elements of $I$, 
and the summation ranges over all subsets $I$ of size $i$). 

\begin{theorem}\label{sl-coinv} 
Assume that $q\in\mc$ is not a root of unity. 
\begin{itemize} 
\item[(i)] 
The subset of $\oslqn$-coinvariants in $\omqn$ with respect to $\alpha$ 
is an $N$-variable commutative polynomial subalgebra of $\omqn$ 
generated by $\si$, $i=1,\ldots,N$. 
\item[(ii)] 
The subset of $\oslqn$-coinvariants in $\oslqn$ with respect to $\alpha$ 
is an $(N-1)$-variable commutative polynomial subalgebra of $\oslqn$ 
generated by $\pi(\si)$, $i=1,\ldots,N-1$. 
\item[(iii)] 
The subset of $\oslqn$-coinvariants in $\omqn$ with respect to $\beta$ 
is an $N$-variable commutative polynomial subalgebra of $\omqn$ 
generated by $\ti$, $i=1,\ldots,N$. 
\item[(iv)] 
The subset of $\oslqn$-coinvariants in $\oslqn$ with respect to $\beta$ 
is an $(N-1)$-variable commutative polynomial subalgebra of $\oslqn$ 
generated by $\pi(\ti)$, $i=1,\ldots,N-1$. 
\end{itemize} 
\end{theorem} 

\begin{proof} Both $\alpha$ and $\beta$ are homogeneous of degree zero, 
so the subsets of $\oslqn$-coinvariants in $\omqn$ 
are the same as for the corresponding $\oglqn$-coaction 
by Proposition~\ref{rel-inv}. 
Hence (i) and (iii) follow from the main result of \cite{dl}, 
describing the spaces of $\oglqn$-coinvariants.  
The statements (ii) and (iv) are immediate consequences of (i) and (iii).  
Indeed, since $\oslqn$ is cosemisimple by our assumption on $q$, 
the subsets of $\oslqn$-coinvariants in $\omqn$ are mapped 
by $\pi$  (which is a morphism of comodules) onto 
the corresponding subsets of coinvariants in $\oslqn$. 
So the only thing left to show is that 
$\pi(\si)$, $i=1,\ldots,N-1$ are algebraically independent, 
and  
$\pi(\ti)$, $i=1,\ldots,N-1$ are algebraically independent. 
This follows from the observation that 
$\ker(\pi)\cap\omqn^{\coa}=\omqn^{\coa}(\detq-1) $ 
and 
$\ker(\pi)\cap\omqn^{\cob}=(\detq-1) \omqn^{\cob}$ 
by Lemma~\ref{proj}. 
\end{proof}


\bigskip 

\noindent M. Domokos: 
R\'enyi Institute of Mathematics, Hungarian Academy of Sciences,\\ 
P.O. Box 127, 1364 Budapest, Hungary\\
E-mail: domokos@renyi.hu, domokos@maths.ed.ac.uk\\
\noindent(Domokos is in Edinburgh until February 2004.)
\\
\\
\noindent T. H. Lenagan: 
School of Mathematics, University of Edinburgh,
\\ James Clerk Maxwell Building, King's Buildings, Mayfield Road, 
\\Edinburgh EH9 3JZ, Scotland
\\E-mail: tom@maths.ed.ac.uk

\end{document}